\documentclass{iopconfser}
\usepackage{graphicx}
\usepackage{amsmath}

\begin{document}
\title{A superconvergence result in the RBF-FD method}

\author{Andrej Kolar-Požun$^{1,2}$, Mitja Jančič$^{1}$ and Gregor Kosec$^{1}$}
\affil{$^1$ Parallel and Distributed Systems Laboratory, Jožef Stefan Institute, Ljubljana, Slovenia}
\affil{$^2$ Faculty of Mathematics and Physics, University of Ljubljana, Ljubljana, Slovenia}

\email{andrej.pozun@ijs.si, mitja.jancic@ijs.si, gregor.kosec@ijs.si}
\begin{abstract}
Radial Basis Function-generated Finite Differences (RBF-FD) is a meshless method that can be used to numerically solve partial differential equations. The solution procedure consists of two steps. First, the differential operator is discretised on given scattered nodes and afterwards, a global sparse matrix is assembled and inverted to obtain an approximate solution. Focusing on Polyharmonic Splines as our Radial Basis Functions (RBFs) of choice, appropriately augmented with monomials, it is well known that the truncation error of the differential operator approximation is determined by the degree of monomial augmentation. Naively, one might think that the solution error will have the same order of convergence. We present a superconvergence result that shows otherwise - for some augmentation degrees, order of convergence is higher than expected.
\end{abstract}

\section{Introduction}
Some of the most popular topics in numerical analysis are arguably numerical solution procedures for Partial Differential Equations (PDEs). Both because this field has numerous applications, as PDEs describe a variety of real-life phenomena, and because of its many challenges - PDEs are typically much more complicated than ordinary differential equations with only very simple cases being solvable analytically.

Without delving too much into the history, let us mention that as far as real-life applications go, the most commonly used methods for numerical solutions of PDEs consider computational domains with some additional structure. For example, in the widespread Finite Element (FEM)~\cite{FEM} and Finite Difference (FDM)~\cite{FDM} Methods, this structure is a (typically triangular) mesh and a regular grid, respectively. 

The method that we focus on in this paper instead falls into the scope of meshless methods that work directly with scattered nodes without any additional structure and with very mild requirements on their geometry. The main advantage of this approach is its greater flexibility (in particular when it comes to adaptive solution procedures~\cite{hp}), as we eliminate a possibly cumbersome step of meshing, while maintaining an accurate description of irregular geometries.

Over the last few decades, several different meshless methods have been proposed. We refer an interested reader to~\cite{review} for a more in-depth survey. We will focus on the Radial Basis Function-generated Finite Difference (RBF-FD) method, proposed in the early 2000s~\cite{rbffd} that is particularly attractive due to its simplicity. It is a strong form method that directly generalises the FDM by approximating a linear differential operator locally on a given neighbourhood of a computational node. As its name implies, it is based on the Radial Basis Functions (RBFs), which have historically been particularly attractive for scattered data interpolation due to their positive definiteness properties~\cite{fasshauer}. 

In our studies, we opt for the Polyharmonic Splines (PHS) as our RBF of choice. Concretely, the radial cubics $\phi(r) = r^3$ will be our RBF basic function. For clarity, we will refer to the employed solution procedure as PHS RBF-FD. Compared to smooth RBFs, PHSs posses a weaker notion of conditional positive definiteness, which requires us to add some number of monomials to the approximation, a procedure known as monomial augmentation. Nevertheless the advantage gained from using such RBFs is that they do not posses a shape parameter, which makes them easier to use. For more information, we once again refer the reader to~\cite{fasshauer}.

Specifically for PHS RBF-FD, various studies of the method's properties have already been performed. In a paper by Flyer et al.~\cite{role1} it is shown that the PHS order is irrelevant as far as the convergence order of the approximation is concerned. In fact, the latter is instead determined by the degree of the augmenting monomials. This result has been further developed by Bayona~\cite{insight}, who analysed monomials' effect in greater detail and derived a formula expressing the PHS RBF-FD truncation error. 

PHS RBF-FD solution procedure is reminiscent of the standard FDM - The linear differential operator appearing in the PDE is discretised locally at each computation point. This allows us to transform the PDE into a set of algebraic equations, which is then solved to obtain a numerical solution. Usually, the order of convergence of such a solution procedure matches the order of operation approximation error.
However, in some special cases, the so called superconvergence can occur, where the solution converges at a higher order than expected.
Superconvergence has been well researched in the context of finite element methods~\cite{superconvergence}, but has also been observed in methods that are variants of the FDM, such as in the Shortley-Weller approximation, where differing truncation orders at the boundary and the interior of the domain cause a superconverge result~\cite{shortley}. More recently, another superconvergence result was presented~\cite{quadrature}, where finite difference schemes, derived from weak form with specially chosen quadrature points cause certain error cancellations leading to the convergence order elevation.

In this paper we start by applying PHS RBF-FD to solve a common test problem - The Poisson equation on a unit disc, albeit irregularly discretised. We first show that the operator truncation error behaves as expected. However, for even augmentation degrees, the solution convergence order is of approximately one higher than expected.

In the following section our problem setup and the methods used are further described.
In section 3 our main result is shown and briefly analysed.
We conclude with section 4, repeating the main remarks and providing possible ideas for further research.
\section{Problem setup}
We consider the Poisson problem:
\begin{equation} \label{eq:PDE}
\nabla^2 u(x,y) = f(x,y),
\end{equation} 
for $x^2+y^2 \leq 1$, where $f(x,y)$ is chosen such that the analytical solution is
\begin{equation} \label{eq:solution}
u(x,y) = 1 + \sin(4x) + \cos(3x) + \sin(2y).
\end{equation}
$u(x,y)$, restricted to $x^2+y^2 = 1$ also gives us the Dirichlet boundary conditions.

\begin{figure}[h]
	\centering
	\includegraphics[width=.5\textwidth]{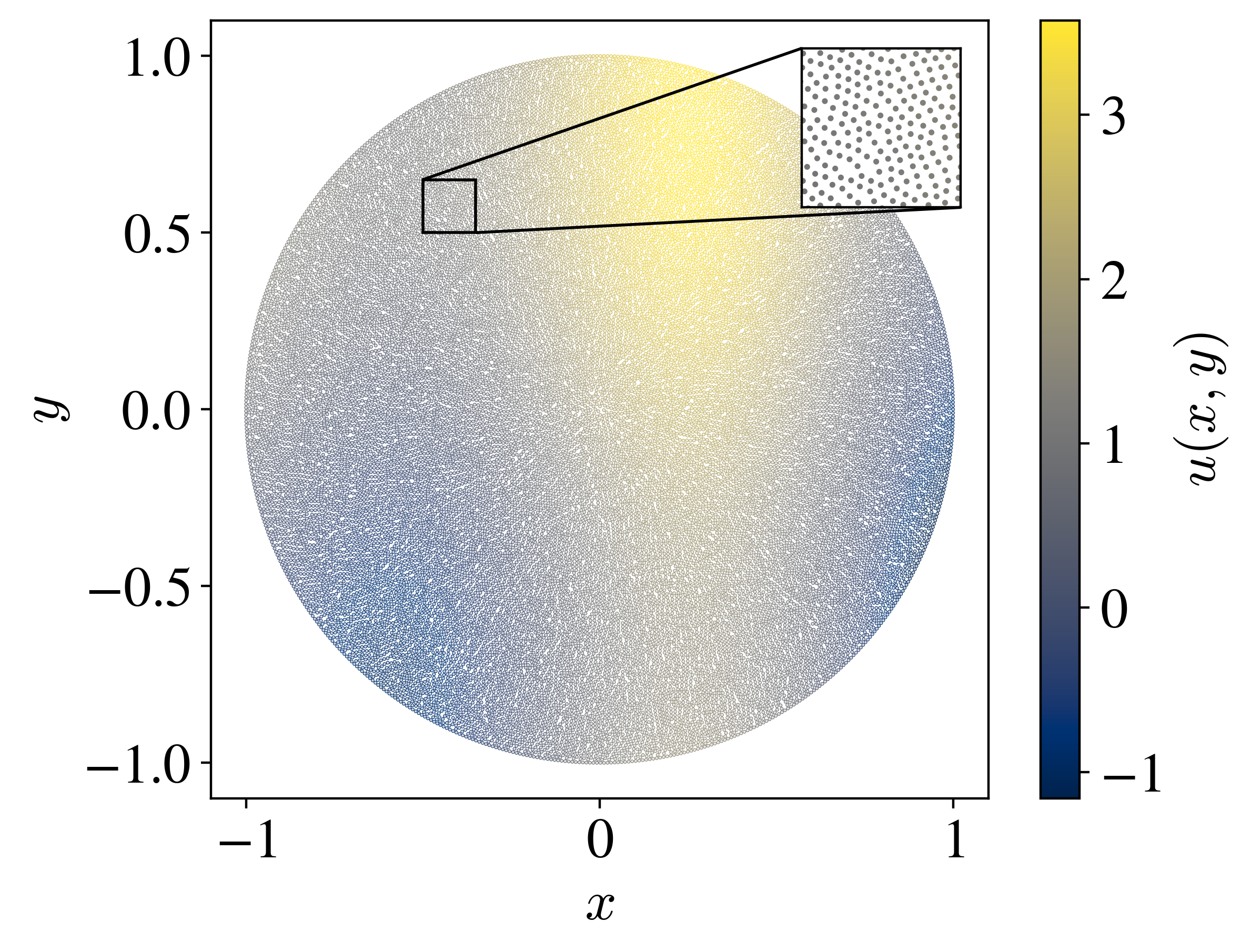}
	\caption{Analytical solution $u(x,y)$, shown on an example discretisation generated during our solution procedure.}
	\label{fig:discretisation}
\end{figure}

As mentioned, the problem will be solved with the PHS RBF-FD method. First, the domain is discretised with the algorithm described in~\cite{fill}, providing a locally regular discretisation, with a desired internodal distance $h$. An example discretisation can be seen on Figure~\ref{fig:discretisation}.
To each discretisation point $(x_i,y_i)$ we then associate the set of its $n$ closest nodes, known as its stencil $S_i$. We set $n = (m+1)(m+2)$, where $m$ is the PHS RBF-FD monomial augmentation degree. This is a common choice and can have a positive impact on the method stability, as explained in~\cite{role2}.

A local approximation of $\nabla^2$ over the stencil is then obtained by the standard RBF-FD procedure and involves solving an appropriate linear system of equations. As this paper deals with a specific property of the RBF-FD method, we will assume the reader already has some familiarity with it and not delve into the details. If that is not the case,~\cite{explainRBF} should fill in the missing details.
We then obtain an approximation of $\nabla^2$ in the following form:
\begin{equation} \label{eq:disc}
\nabla^2 u(x_i,y_i) \approx \sum_{j \in S_i} w_j u(x_j,y_j),
\end{equation}
where $w_j$ are the RBF-FD weights, expressing $\nabla^2$ as a linear combination of function values in the the stencil $S_i$. Writing Equation~\ref{eq:disc} for each node in the domain interior, and an equation expressing the Dirichlet boundary conditions for each node on the boundary, we obtain a global (and, due to the locality of the approximation, sparse) system, which we then solve to obtain a numerical solution. The BiCGSTAB iterative solver has been used for that matter. All of the source code can be found on our public git repository\footnote{https://gitlab.com/e62Lab/public/2024\_cp\_eurotherm\_superconvergence}. Note that it relies heavily on our in-house (open source) C++ library Medusa~\cite{medusa}.

In what follows, we will use the just described procedure to numerically solve our model problem and study the error convergence. 
\section{The main result}
Two different errors will be considered. First, the operator approximation error: For each point we compute the difference between the analytical and RBF-FD $\nabla^2$ applied to $u(x,y)$ (Equation~\ref{eq:solution}). For convenience we introduce the vector $\textbf{e}_\mathrm{op}$, its $i$-th component being the error at $(x_i,y_i)$. We will consider the mean error, which by a slight abuse of notation will be denoted by $||\textbf{e}_\mathrm{op}||_1$ (despite it actually being the $\ell_1$ norm, divided by the number of discretisation nodes $N$). All of the coming results have been observed also in $\ell_\infty$ error norm, but will not be displayed for brevity reasons.

\begin{figure}[h]
	\centering
	\includegraphics[width=.6\textwidth]{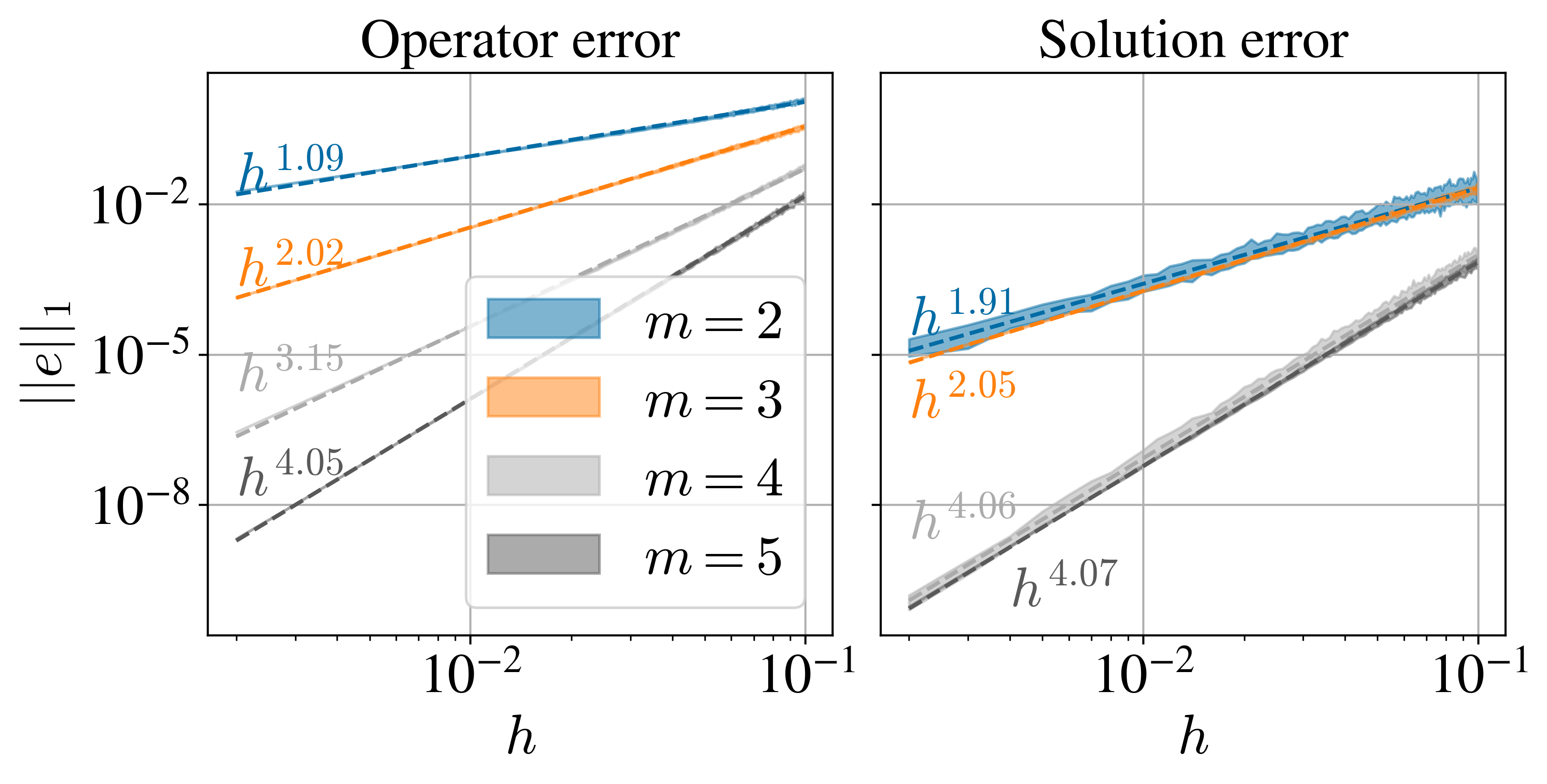}
	\caption{Convergence of mean errors under decreasing $h$. Operator and solution errors are shown on the left and right, respectively for different monomial augmentation degrees $m$. For each $m$, the errors have been computed for 30 different discretisation sets. The results lie in the corresponding shaded regions.}
	\label{fig:convergenceAll}
\end{figure}

The error convergence of the $\nabla^2$ RBF-FD approximation is shown on the left side of Figure~\ref{fig:convergenceAll}.
As our discretisation is irregular, we have computed these errors over 30 different discretisation sets, varied by changing the random seed in our previously mentioned discretisation algorithm.
We can see that for each monomial augmentation degree $m$, the errors scale as approximately $\propto h^{m-1}$, which is in agreement with the theory - operator discretisation error should scale as $ \propto h^{m+1-l}$, where $l$ is the derivative order, in our case $l=2$~\cite{explainRBF}. 
Next we consider the solution error, expressing the error of the whole solution procedure when solving our chosen Poisson problem. Analogously to the previous notation we introduce the vector $\textbf{e}_\mathrm{sol}$. Solution error convergence can be seen on the right side of Figure~\ref{fig:convergenceAll}, which shows our main observation of the paper -  convergence orders match those of the operator approximation for odd augmentation degrees $m$, where the error scales as $\propto h^{m-1}$, but we get an order higher for even $m$ - approximately $\propto h^m$. Additionally, the error spread for even $m$ is noticeably larger, hinting at a different underlying mechanism.

\begin{figure}[h]
	\centering
	\includegraphics[width=.6\textwidth]{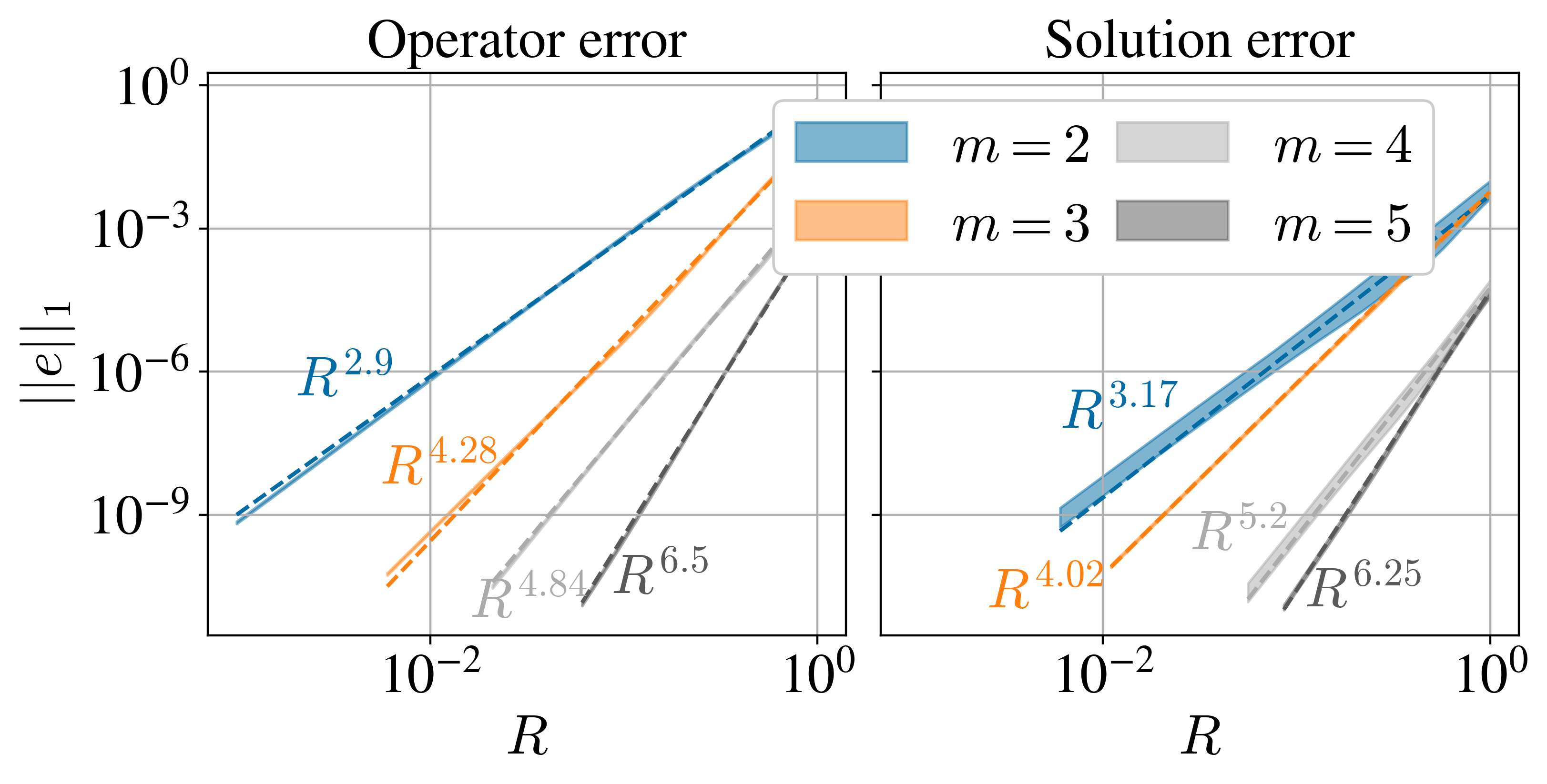}
	\caption{Solution and operator error dependences on scaling factor $R$ shown for different monomial augmentation degrees $m$. For each $m$, the errors have been computed for 30 different discretisation sets. The results lie in the corresponding shaded regions.}
	\label{fig:convergenceR}
\end{figure}

For our analysis, we first consider an alternative procedure of studying error convergence. Keeping $h$ fixed at $0.05$, we instead make a transformation $u(x,y) \to u(Rx,Ry)$ and study the error scaling as we increase $R$, effectively flattening the chosen solution. Error convergence with respect to $R$ can seen on Figure~\ref{fig:convergenceR}. Both the operator and solution errors now scale as $R^{m+1}$, which is expected.

The lack of order elevation when studying the error dependence on $R$ implies that for even augmentation degrees $m$ we do not actually gain any extra approximation power (as in, method suddenly becoming exact for higher degree monomials). Instead, an extra factor of $h$ must appear somewhere in the final step of the solution procedure. Since that step involves solving a sparse system with the global matrix $A$, it makes sense to look at it more closely.


We now attempt to study the error convergence more systematically by applying Bayona's explicit formula~\cite{insight} for the RBF interpolation error over a stencil around the point $(x_0,y_0)$:
\begin{equation}
\begin{split}
\hat{u}(x,y) - u(x,y) = \sum_{k=s+1}^\infty &L_k [u(x_0,y_0)] \big(\textbf{p}_k^T \boldsymbol{\hat{\psi}}(x,y) - p_k(x-x_0,y-y_0) - \\ &-\textbf{p}_k^T W (P^T \boldsymbol{\hat{\psi}}(x,y) - \textbf{p}(x-x_0,y-y_0)) \big),
\end{split}
\end{equation}
where the notation is similar to Bayona's and we do not repeat it here.

\begin{figure}[h]
	\centering
	\includegraphics[width=.6\textwidth]{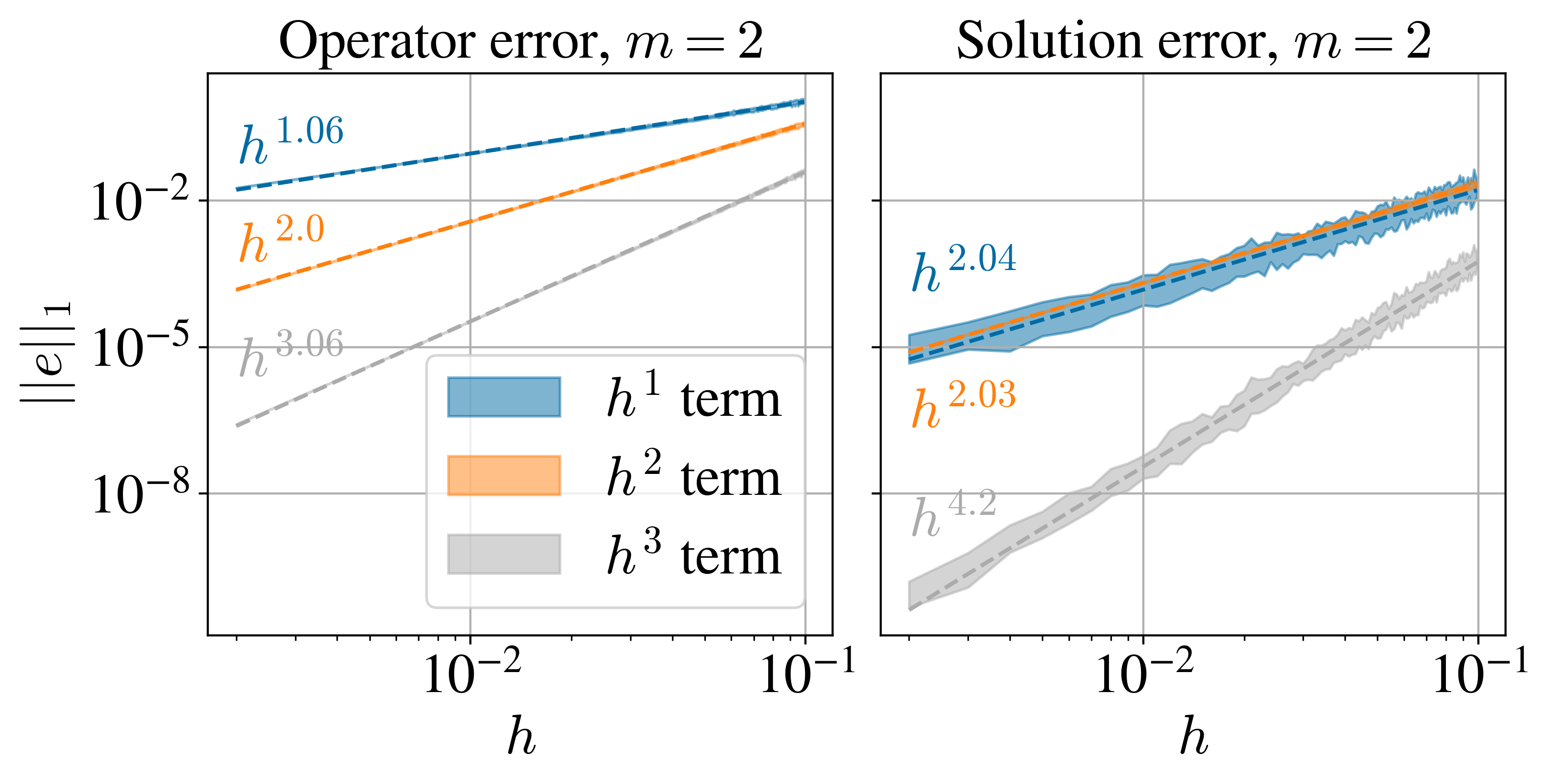}
	\caption{Average error for each term in Bayona's formula for the operator approximation on the left, and the same terms after solving the global system on the right. Once again, the results over 30 different discretisations are shown. In both cases $m=2$ was fixed and first few leading terms in the formula were considered.}
	\label{fig:errorFormula}
\end{figure}

Applying $\nabla^2$ over the above and evaluating at some discretisation point $(x_i,y_i)$, we obtain a simplified formula for the operator approximation error:
\begin{equation} \label{eq:bayona}
\nabla^2\hat{u}(x_i,y_i) - \nabla^2u(x_i,y_i) = \sum_{k=s+1}^\infty L_k [u(x_i,y_i)] (\textbf{p}_k^T \cdot \textbf{w}),
\end{equation}
where we have expressed the vector of $\nabla^2$ RBF-FD weights $\textbf{w}$ as $\textbf{w} = \nabla^2 \boldsymbol{\psi}(x,y)$ and used the facts that $\boldsymbol{\psi}(x,y) = (I - W P^T) \boldsymbol{\hat{\psi}}(x,y) + W \textbf{p}(x,y)$ and $p_k(0,0) = 0$~\cite{insight}. The weights $\textbf{w}$ scale as $h^{-2}$, while $\textbf{p}_k^T$ is a vector of monomials scaling as, for example $h^{m+1}$ for the first terms in the error, and with higher exponents for later terms. We can use this formula to analyse the operator error term by term, grouped by powers of $h$. The results are seen on the left side of Figure~\ref{fig:errorFormula} and match the previously observed operator error scalings.

Next, we use this formula to also analyse the behaviour of the solution error. Consider the RBF-FD approximation of $\nabla^2$, applied to the solution error: $\nabla^2 \textbf{e}_\mathrm{sol}(i) = \nabla^2 (\hat{u}(\textbf{x}_i) - u(\textbf{x}_i)) = f_i - \nabla^2 u(\textbf{x}_i) = \textbf{e}_\mathrm{op} (i)$. This, together with the fact that the solution error on the boundary is equal to $0$ (due to the Dirichlet boundary conditions) implies
\begin{equation} \label{eq:errors}
A \textbf{e}_\mathrm{sol} = \textbf{e}_\mathrm{op}.
\end{equation}
Therefore, we can analyse how each of the terms in Bayona's error formula behave also after the global system is inverted - we simply plug the corresponding term as the RHS. Right side of Figure~\ref{fig:errorFormula} shows that terms, that originally contained an odd power of $h$, obtain an extra factor of $h$ after the system is inverted, while the ones with even powers are unaffected. This is why we observe superconvergence only for even augmentation degrees, as in those cases the dominant error term contains odd powers of $h$, while for odd augmentation degrees, the opposite applies. 

Summarising, we have demonstrated that the observed order elevation can already be seen by considering Bayona's error formula, as written in Equation~\ref{eq:bayona}. Further study of the latter could be easier than of the original problem (Equation~\ref{eq:PDE}), as we have essentially reduced it to the study of the solution to a linear system, where the functional forms of the matrix $A$ and the RHS are known (Equation~\ref{eq:errors}). Further investigation is difficult, as it would require an understanding of what precisely happens during the process of inverting the matrix $A$, which is a challenge due to its size. Some ideas are mentioned in the conclusion below.

\section{Conclusion}
This paper considered the error analysis of the RBF-FD method, applied to the Poisson problem on a unit disc.
We have opted for radial cubics as our RBF with monomial augmentation, for which it is well established on how the operator error scales with the internodal distance $h$. 

Despite the simplicity of the model problem, we have noticed unexpected behaviour of the solution error convergence: for even augmentation degrees, the convergence order is of approximately one higher than predicted by the theory.

We have managed to simplify the problem, showing that we can use Bayona's error formula to study the error behaviour systematically term by term.
 
There are multiple ways of proceeding in this research that we are actively attempting in our further studies. First, one can make use of the known theory mentioned in the introduction, which other authors studying similar phenomena have successfuly applied. More tests of purely experimental nature are also being performed to verify which properties of the problem or its solution procedure affect the observed orders of convergence.

\section*{Acknowledgments}
The authors acknowledge the financial support from the Slovenian Research and Innovation Agency (ARIS) in the framework of the research core funding No. P2-0095 and the Young Researcher program PR-12347.

\bibliographystyle{abbrv}
\bibliography{Manuscript}
\end{document}